\input amstex
\documentstyle{amsppt}
\refstyle{C}
\TagsOnRight
\define\e{@,@,@,@,@,}
\define\img{{\roman{Im}}\,}
\define\g{\varGamma}
\define\gv{\varGamma_v}
\define\ng{\negthickspace}
\define\x{\delta}
\define\y{\delta^{@,@,\prime}}
\define\xv{\delta_v}
\define\yv{\delta_v^{@,@,\prime}}
\define\pic{\roman{Pic}}
\define\iv{\roman{Div}}
\define\br{\roman{Br}}
\define\kb{\bar{K}}
\define\cok{\roman{Coker}\,}
\define\ra{\rightarrow}
\define\krn{\roman{Ker}\,}
\define\dgr{\roman{deg}}
\define\lb{\lbrack @,@,@,}
\define\rb{@,@,@,\rbrack}
\define\be{@!@!@!@!@!}
\define\te{\Delta}
\define\tprod{\tsize{\prod}}
\define\ue{\Delta^{\be\prime}}
\font\tencyr=wncyr10
\font\sevencyr=wncyr7
\font\fivecyr=wncyr5
\newfam\cyrfam
\def\cyr{\fam\cyrfam\tencyr}
\textfont\cyrfam=\tencyr
\scriptfont\cyrfam=\sevencyr
\scriptscriptfont\cyrfam=\fivecyr
\loadbold

\topmatter
\title Brauer groups and Tate-Shafarevich groups\endtitle
\author Cristian D. Gonzalez-Aviles\endauthor
\affil Facultad de Ciencias, Universidad de Chile\endaffil
\address{Casilla 653, Santiago, Chile.
\email cgonzale\@uchile.cl\endemail}\endaddress
\bigskip

\thanks{Supported by Fondecyt through grants 1981175 and
1000814.}\endthanks

\abstract{Let $X_K$ be a proper, smooth and geometrically connected curve
over a global field $K$. In this paper we generalize a formula of Milne
relating the order of the Tate-Shafarevich group of the Jacobian of $X_K$
to the order of the Brauer group of a proper regular model of $X_K$.
We thereby partially answer a question of Grothendieck.}
\endabstract

\subjclass 11G35, 14K15, 14G25\endsubjclass
\keywords{Brauer groups, Tate-Shafarevich groups, Jacobian variety,
index of a curve, period of a curve}\endkeywords
\endtopmatter

\nologo

\heading 0. Introduction\endheading

Let $K$ be a global field, i.e. $K$ is a finite extension of
${\Bbb Q}$ (the ``number field case") or is finitely generated and of
transcendence degree 1 over a finite field (the ``function field case").
In the number field case we let $U$ denote a nonempty open subscheme of
the spectrum of the ring of integers of $K$, and when $K$ is a function
field in one variable with finite field of constants $k$, we
let $U$ denote a nonempty open subscheme of the unique smooth complete
curve over $k$ whose function field is $K$. We will write $S$ for the
set of primes of $K$ not corresponding to a point of $U$. Further, we
will write $\kb$ for the separable algebraic closure of $K$ and $\g$ for
the Galois group of $\kb$ over $K$. The completion of $K$ at a prime
$v$ will be denoted by $K_v$.

Assume now that a connected, regular, 2-dimensional scheme $X$ is given
together with a proper morphism $\pi\: X\rightarrow U$ whose generic fiber
$X_K=X\times_U{\roman{Spec}}\, K$ is a smooth geometrically connected curve
over $K$. Let $\x$ (resp. $\y$) denote the {\it{index}} (resp.
{\it{period}}\,) of $X_K$. These integers may be defined
as the least
positive degree of a divisor on $X_K$ and the least positive
degree of a divisor class in $\pic(X_{\kb})^{\g}$, respectively,
where $X_{\kb}=X_K\otimes_K \kb$. For each prime
$v$ of $K$, we will write $\xv$ and $\yv$ for the analogous
quantities associated to the curve $X_{K_v}$. It is known that there
are only finitely
many primes $v$ for which $\xv\neq 1$. Further, Lichtenbaum [21] has
shown that $\xv$ equals either $\yv$ or $2\e\yv$ for
each prime $v$. We will write $d$ for the number of primes $v$
for which $\xv=2\e\yv$. These primes were called ``deficient" in
[35], and we sometimes refer to $d$ as {\it the number of deficient
primes of $X_K$}. Now let $A$ be the Jacobian variety of $X_K$. It has
long been known (see [46], \S 3) that there exist close
connections between
the Brauer group ${\br}(X)$ of $X$ and the Tate-Shafarevich
group ${\cyr W}(A)$ of $A$. These connections were explored at length
by Grothendieck in his paper [15], in a more general setting than the
one considered here. The following theorem can be extracted from [15],
p. 121.

\proclaim\nofrills{Theorem\, }{\rm (Grothendieck)}. In the function
field case, suppose that $\xv=1$ for all primes $v$ of $K$. Then
there exist a finite group $T_2\subset{\cyr W}(A)$ of order $\y$,
a finite group $T_3$ of order dividing $\x$, and an exact sequence
$$0\rightarrow{\br}(X)\rightarrow{\cyr W}(A)/T_2\rightarrow
T_3\rightarrow 0.$$
\endproclaim

Regarding this result, Grothendieck (op. cit., p. 122) asked for the
exact order of $T_3$. This problem was solved by Milne in [28]
(see also [29], III.9.6), who used Cassels-Tate duality to compute
the order of $T_3$. In order to state Milne's result, which covers
both the function field and number field cases, we need the following
definition. Let
$${\br}(X)'={\krn}\!\Biggl[@,@,\br(X)\ra\bigoplus_{v\in S}
\br(X_{K_v})@,@,\Biggr ].$$
Thus ${\br}(X)'$ is the group of elements in the Brauer group
of $X$ becoming trivial on $X_{K_v}$ for all $v\in S$. Then one has the

\proclaim\nofrills{Theorem\, }{\rm (Milne)}. Assume that
$\xv=1$ for all
primes $v$ of $K$ and that ${\cyr W}(A)$ contains no nonzero
infinitely divisible elements. Then the period of $X_K$ equals its
index, i.e. $\y=\x$, and there exist finite groups $T_2$ and
$T_3$ of order $\delta$ and an exact sequence
$$0\rightarrow{\br}(X)'\rightarrow{\cyr W}(A)/T_2\rightarrow
T_3\rightarrow 0.$$
In particular, if one of ${\cyr W}(A)$ or ${\br}(X)'$ is
finite, then so is the other, and their orders are related by
$$\x^{\e 2}\e\lb\br(X)'\rb =\lb{\cyr W}(A)\rb.$$
\endproclaim

Grothendieck (loc. cit.) went on to pose the problem of making
explicit the relations between $\br(X)$ and ${\cyr W}(A)$
when the integers $\xv$ are no longer assumed to be equal to
one\footnote{To be precise, Grothendieck stated an equivalent
form of this problem.}. In this paper we generalize the methods developed
by Milne in [28] to prove the above theorem and obtain the following
stronger result, which may be viewed as a partial solution
to Grothendieck's problem
\footnote{In the last section of the paper we show that if
a certain plausible conjecture holds, then we can give a complete
answer to Grothendieck's question.}.

\proclaim{Main Theorem} Assume that the integers $\yv$ are
relatively prime in pairs (i.e., $(\yv,@,@,@, \delta^{@,@,\prime}_
{\bar v})=1$ for all $v\neq{\bar v}$) and that ${\cyr W}(A)$ contains no
nonzero infinitely divisible elements. Then there is an exact sequence
$$0\rightarrow T_{0}\rightarrow T_1\rightarrow\br(X)'
\rightarrow{\cyr W}(A)/T_2\rightarrow T_3\rightarrow 0\,,$$
in which $T_0$, $T_1$, $T_2$ and $T_3$ are finite groups of orders
$$\align \lb T_{0}\rb &=\delta/\delta'\\
\lb T_{1}\rb &=2^e \\
\lb T_{2}\rb &=\y/{\tsize{\prod}}\e\yv\\
\lb T_{3}\rb &=\dfrac{\y/\prod\be\yv}{2^f}\,,\endalign$$
where
$$e={\roman{max}}(0, d-1)$$
and
$$f=\cases 1\qquad\text{if $\;\y/\prod\be\yv$ is even and $d\geq 1$}\\
0\qquad\text{otherwise.}\endcases$$
Here $d$ is the number of deficient primes of $X_K$ defined previously.
In particular, if one of ${\cyr W}(A)$ or $\br(X)'$ is
finite, then so is the other, and their orders are related by
$$\x\y\,\lb\br(X)'\rb =2^{@,@,@, e+f}\,{\tsize{\prod}}
(\yv)^2\; \lb{\cyr W}(A)\rb .$$
\endproclaim
Some immediate corollaries are
\proclaim{Corollary 1} Suppose that one of ${\cyr W}(A)$ or
$\br(X)'$ is finite. Assume also that $\xv=\yv$ for
all $v$ (which holds for instance if $X_K$ has genus 1
\footnote{Regarding the parenthetical remarks in the statements of
both corollaries, see [21]}), and that
these integers are relatively prime in pairs. Then $\x=\y$,
and
$$\x^{\e 2}\,\lb\br(X)'\rb={\tsize{\prod}}\,\xv^{\e 2}\;
\lb{\cyr W}(A)\rb.$$
In particular if $\xv=1$ for all $v$, then
$$\x^{\e 2}\,\lb\br(X)'\rb =\lb{\cyr W}(A)\rb.$$
\endproclaim

Note that the last formula in the statement of Corollary 1 is
precisely the formula of Milne stated before.

\proclaim{Corollary 2} Assume that one of ${\cyr W}(A)$ or $\br(X)'$
is finite and that $\yv=1$ for all $v$ (the latter holds
for instance if $X_K$ has genus 2). Then
$$\x\y\,\lb\br(X)'\rb =2^{@,@,@, e+f}\,\lb{\cyr W}(A)\rb,$$
where $e$ and $f$ are as in the statement of the theorem.
\endproclaim

In the function field case our main result, when combined with a formula
of Gordon [10, p.$\hphantom{.}196$], should imply the expected equivalence of the
conjectures of Artin and Tate for $X$ [46, Conj. C] and Birch and
Swinnerton-Dyer for $A$ [46, Conj. B] (this is Tate's ``elementary"
conjecture (d) of [46]), at least under the additional assumption that
the structural morphism $\pi\: X\rightarrow U$ is cohomologically flat in
dimension 0. We expect to address this issue in a separate publication.
\bigskip
\heading Acknowledgements \endheading

I am grateful to the people of Fondecyt for their financial
support and their patience. I also thank the library staff
at IMPA (Rio de Janeiro, Brasil) for their bibliographical
assistance.

\heading 1. Preliminaries\endheading

We keep the notations introduced in the previous section. Thus
in particular $X$ is a regular connected scheme of dimension 2
equipped with a proper morphism $\pi\: X\rightarrow U$ whose generic
fiber $X_K$ is a smooth geometrically connected curve over $K$.

\smallskip

{\it{Remark}}. The $U$-scheme $X$ is a proper regular model of
its generic fiber. Conversely, if we start with a (geometrically
connected) proper and smooth curve $X_K$ over $K$, then there is a
closed immersion $X_K\rightarrow {\Bbb{P}}^{@,@,@,n}_{\ng K}$
for some $n$ and we can obtain a $U$-scheme $X$ as above, with generic
fiber $X_K$, by applying Lipman's desingularization process [2], [22] to
the schematic image of $X_K$ in ${\Bbb P}^{@,@, n}_{\! U}$.

\smallskip

The Picard scheme of $X_K/K$, $\pic_{X_{\be K}\be/\be K}$, is a
smooth group
scheme over $K$ whose identity component, $\pic^0_{X_{\be K}\be/
\be K}$, is
an abelian variety, the Jacobian variety of $X_K$. Henceforth,
we will write $P$ for $\pic_{X_K/K}$ and (in accordance with
previous notations) $A$ for $\pic^0_{X_{\be K}\be/\be K}$. The following
holds. If $L$ is any field containing $K$ such that $X_K(L)$ is
nonempty, then $P(L)=\pic(X_L)$ (for the basic facts on the relative
Picard functor, see [4, \S 8.1] or [12]). Further, there is an exact
sequence of $\g$-modules

$$0\ra A(\kb)\ra P(\kb)\overset\dgr\to{\longrightarrow}\Bbb Z\ra 0,$$
where $\dgr$ is the degree map on $P(\kb)=\pic(X_{\kb})$. In
particular $A(\kb)$ may be identified with $\pic^0(X_{\kb})$, the
subgroup of $\pic(X_{\kb})$ consisting of divisor classes of degree zero
(in this paper we shall regard the elements of $\pic(X_{\kb})$ mainly as
classes of divisors. See below). Now we observe that $P(K)=P(\kb)^\g=
\pic(X_{\kb})^\g$. Further, there is an exact sequence
(deduced from the preceding one by taking $\g$-invariants)

$$0\ra A(K)\ra P(K)\overset\dgr\to{\longrightarrow}\y\e\Bbb Z\ra 0\,,$$
where $\y$ is the period of $X_K$ as defined previously.

For any regular connected scheme $Y$, $R(Y)$ will denote the field
of rational functions on $Y$. There is an exact sequence
$$0\ra R(X_{\kb})^*/\kb^*\ra\iv^0(X_{\kb})\ra\pic^0(X_{\kb})\ra 0$$
which induces an exact sequence

$$0\ra R(X_K)^*/K^*\ra\iv^0(X_K)\ra A(K)\ra A(K)/\pic^0(X_K)\ra 0.\tag 1$$
Note that $A(K)/\pic^0(X_K)$ is a finite abelian group since
it is finitely generated (by the Mordell-Weil theorem) and isomorphic
to a subgroup of the torsion group $H^1(\g, R(X_{\kb})^*/\kb^*)$.
Similarly, $P(K)/\pic(X_K)$ is finite and there is an exact sequence

$$0\ra R(X_K)^*/K^*\ra\iv(X_K)\ra P(K)\ra P(K)/\pic(X_K)\ra 0.\tag 2$$
\newpage

We also have an exact sequence

$$0\ra \iv^0(X_K)\ra\iv(X_K)\overset\dgr\to{\longrightarrow}\x\,
\Bbb Z\ra 0\,,$$
where $\x$ is the index of $X_K$ as defined previously. We note
that $\x$ may also be defined as the greatest common
divisor of the degrees of the fields $L$ over $K$ such that
$X_K(L)\neq\emptyset$.

\proclaim{Lemma 1.1} The period of $X_K$ divides its index, i.e.
$\y|\,\x$, and
$$\lb P(K):\pic(X_K)\rb=(\x/\y)\cdot\lb A(K):\pic^0(X_K)\rb.$$
\endproclaim
{\it Proof}. The first assertion of the lemma follows from the
definitions. Regarding the second, an application of the
snake lemma to the commutative diagram
$$\minCDarrowwidth{.5cm}
\CD 0 @>>>\iv^0(X_K) @>>>\iv(X_K) @>>>\x\,\Bbb Z @>>> 0\\
      @.     @VVV           @VVV        @VVV            @.\\
      0 @>>> A(K)      @>>> P(K)     @>>>\y\,\Bbb Z @>>> 0\endCD $$
yields, using (1) and (2) above, an exact sequence
$$0\ra A(K)/\pic^0(X_K)\ra P(K)/\pic(X_K)\ra \y\,\Bbb Z/\x\e\Bbb Z
\ra 0.$$
The lemma is now immediate. \qed

It is clear from the definitions that if $L$ is any field
containing $K$, then the index (resp. period) of $X_L$ divides the index
(resp. period) of $X_K$. Thus for any prime $v$ of $K$, $\xv|\,\x\,$
and $\,\yv|\,\y$, where $\xv$ and $\yv$ are, respectively, the index and
period of $X_{K_v}$. Further, since $X_K(K_v)\neq\emptyset$ for
all but finitely many primes $v$ (see for example [18, p.$\hphantom{.}
249$, Remark 1.6\,]), we conclude that there are only finitely many
primes $v$ for which $\xv\neq 1$. The analogous statement with $\xv$
replaced by $\yv$ holds true as well, since $\,\yv|\,\xv$ for each
$v$. We now write $\te$ (resp. $\ue$) for the least common
multiple of the integers $\xv$ (resp. $\yv$). Clearly,
$\te|\,\x$ and $\ue|\,\y$. Further, since $\xv=\yv\;\,\text{or}\;\;
2\e\yv\,$ for each $v$ as mentioned earlier, we have

$$\frac{\te}{\ue}=\cases 1\qquad\text{if $\xv=\yv$ for all $v$}\\
2\qquad\text{otherwise.}\endcases$$
\medskip

We now consider the map
$$\Sigma\, \:\bigoplus_v\e \xv^{-1}\e\Bbb Z/\Bbb Z\ra\Bbb Q/\Bbb Z\,,
\quad (x_v)\mapsto\sum x_v\,.$$

\proclaim{Lemma 1.2} We have
$$\lb\krn\Sigma\rb =\tprod\e\xv/\te$$
and
$$\img\Sigma=\te^{\be -1}\e\Bbb Z\,/\,\Bbb Z.$$
\endproclaim
{\it Proof}. That the image of $\Sigma$ is exactly $\te^{\be -1}\e\Bbb Z
\e/\e\Bbb Z$ follows easily from the fact that ${\roman{gcd}}
\fracwithdelims(){\te}{\xv}=1$. The rest of the lemma is clear.\qed

Next we consider the map
$$D\,\: \Bbb Z/\y\e\Bbb Z\ra\bigoplus_v\e\Bbb Z/\yv\e\Bbb Z$$
given by $\, x\ng\mod{\y}\mapsto(x\ng\mod{\yv})$.
Since the kernel of $D$ is $\ue\e\Bbb Z/\y\e\Bbb Z$, the following
lemma is clear

\proclaim{Lemma 1.3} We have
$$\lb\krn D\rb = \y/\ue$$
and
$$\hphantom{1}\lb\cok D\rb =\tprod\e\yv/\ue.$$
\endproclaim
\medskip

For each prime $v$ of $K$ we will write $\gv$ for the Galois group of
$\kb_v$ over $K_v$.

From the cohomology sequence associated to the
exact sequence of $\gv$-modules $0\ra A(\kb_v)\ra P(\kb_v)\overset\dgr
\to{\longrightarrow}\Bbb Z\ra 0$ we get an exact sequence

$$0\ra\Bbb Z/\yv\e\Bbb Z\ra H^1(\gv,A)\ra H^1(\gv,P)\ra 0.\tag 3$$

Now by work of Lichtenbaum [21] and Milne [26]
(see also [28, Remark I.3.7]), there exists a perfect pairing

$$H^0(\gv,A)\times H^1(\gv,A)\ra\Bbb Q/\Bbb Z\,,$$
where $H^0(\gv,A)$ denotes $A(K_v)/N_{\be\kb_v\be/\be
K_v}\! A(\kb_v)$ if $v$ is archimedean and $A(K_v)$ otherwise.
Relative to this pairing, the annihilator of the image of
$\pic^0(X_{K_v})$ in $H^0(\gv,A)$ under the canonical map $A(K_v)\ra
H^0(\gv,A)$ is exactly the image of $\Bbb Z/\yv\e\Bbb Z$ in $H^1(\gv,A)$
under the map in (3). Consequently the following holds

\proclaim{Lemma 1.4} We have
$$\lb A(K_v):\pic^0(X_{K_v})\rb =\yv.$$
\endproclaim

We now combine the previous lemma with an analogue of Lemma 1.1 to
obtain

\proclaim{Lemma 1.5} We have
$$\lb P(K_v):\pic(X_{K_v})\rb =\xv.$$
\endproclaim
\smallskip

{\it Remark}. The archimedean case of Lemma 1.5 was originally
established by Witt in 1935 [50]. For an interesting review of this and
other related classical results in terms of \'etale cohomology, see
[38, \S 20.1].
\smallskip

We now derive a slight variant of the snake lemma (Proposition 1.6
below). It is one of the basic ingredients of the proof of the Main
Theorem.

Consider the following exact commutative diagram in the category of
abelian groups

$$\minCDarrowwidth{.6cm}
\CD 0 @>>> A_1 @>>> A_2 @>f>> A_3 @>>> A_4 @>>> A_5 @>>> 0\\
      @.   @VVV     @VVV     @VV\eta V     @VV\lambda V @VV\mu V
      @.\\
      0 @>>> B_1 @>>> B_2 @>g>> B_3 @>>> B_4 @>>> B_5 @>>>
0\endCD\tag 4 $$
(we have labeled only those maps which are relevant to our purposes).
We have an induced exact commutative diagram
     
$$\CD  @. A_3 @>>>A_4 @>>>A_5 @>>> 0\\
      @.     @VV{\bar{\eta}}V   @VV\lambda V    @VV\mu V   @.\\
      0 @>>>B_3/\img g @>>> B_4  @>>>B_5 @>>> 0\endCD $$
where ${\bar{\eta}}\,$ is the composition of $\eta$ with the
canonical map $B_3\ra B_3/\img g$. An application of the snake lemma to
the above diagram using the fact that $\img f\subset\krn\bar
{\eta}\,$ yields

\proclaim{Proposition 1.6} To any exact commutative diagram of the form
(4) there is associated an exact sequence
$$0\ra\img f\ra\krn{\bar{\eta}}\ra\krn\lambda\ra\krn\mu\ra
\cok{\bar{\eta}}\,,$$
where ${\bar{\eta}}$ is as defined above.
\endproclaim

The foolowing result supplements Proposition 1.6.

\proclaim{Lemma 1.7} With the above notations, there is an exact
sequence
$$0\ra\krn\eta\ra\krn{\bar{\eta}}\ra\img g\ra\cok\eta\ra\cok{
\bar{\eta}}\ra 0.$$
\endproclaim
{\it Proof}. This is nothing more than the kernel-cokernel sequence
[.., I.0.24] for the pair of maps $A_3\overset\eta\to
{\longrightarrow}B_3\ra B_3/\img g$. The maps in the exact sequence of
the lemma are the natural
ones, e.g. $\img g\ra\cok\eta$ is the composite $\img g
\hookrightarrow B_3\ra B_3/\img\eta=\cok\eta$.\qed

\heading 2. Proof of the Main Theorem\endheading
All cohomology groups below will be either Galois cohomology groups or
\'etale cohomology groups. We will view $\gv$ as a subgroup of $\g$ in
the standard way, i.e., by identifying it with the decomposition group of
some fixed prime of $\kb$ lying above $v$. For each $v$,
$\roman{inv}_v\,\:\br(K_v)\ra\Bbb Q/\Bbb Z$ will denote the usual
invariant map of local class field theory. The $n$-torsion subgroup of
an abelian group $M$ will be denoted by $M_{n-\roman{tor}}$.

We begin by recalling a fundamental exact sequence.
Since $H^q(X_{\kb}, \Bbb G_m)=0$ for all $q\geq 2$ [27, Ex. 2.23(b),
p.$\hphantom{.}$110], the Hochschild-Serre spectral sequence
$$H^p(\g,H^q(X_{\kb},\Bbb G_m))\Rightarrow H^{p+q}(X_K,\Bbb G_m)$$
yields (see [7, XV.5.11]) an exact sequence

$$0\ra\pic(X_K)\ra P(K)\ra\br(K)\ra\br(X_K)\ra H^1(\g,P)\ra 0\,,\tag 5$$
where the zero at the right-hand end comes from the fact that
$H^3(\g,\kb^*)=0$ [29, I.4.21]. (We have used here the well-known
facts that $\pic(X_K)=H^1(X_K,\Bbb G_m)$ and that the Brauer group of
a regular scheme of dimension $\leq 2$ agrees with the cohomological
Brauer group of the scheme.) Similarly, for each prime $v$ of
$K$ there is an exact sequence
$$0\ra\pic(X_{K_v})\ra P(K_v)\overset g_v\to{\longrightarrow}
\br(K_v)\ra\br(X_{K_v})\ra H^1(\gv,P)\ra 0.\tag 6$$

\proclaim{Lemma 2.1} For each prime $v$ of $K$, we have
$$\img g_v = \br(K_v)_{\xv\be -\roman{tor}}.$$
\endproclaim
{\it Proof}. By Lemma 1.5, $\img g_v$ is a subgroup of $\br(K_v)$
of order $\xv$. On the other hand the invariant map
$\roman{inv}_v$ induces an isomorphism $\br(K_v)_{\xv\be -\roman{tor}}
\overset\sim\to{\longrightarrow}\xv^{-1}\e\Bbb Z/\Bbb Z$, whence the lemma
follows.\qed

We now consider the commutative diagram
$$\minCDarrowwidth{.39cm}
\CD \e 0 @>>>\pic(X_K) @>>> P(K) @>>>\br(K) @>>>\br(X_K) @>>> H^1(\g,P)
@>>>\! 0\\
      @.   @VVV     @VVV     @VV\eta V     @VVV  @VVV      @.\\
      0 \! @>>>\ng\underset v\to{\bigoplus}\e\pic(\be X_{K_v}\be)\!
@>>>\ng\underset v\to{\bigoplus}\e P(\be K_v\be)\! @>g>>\ng\underset v\to{
\bigoplus}\e\br(\be K_v\be)\! @>>>\ng\underset v\to{\bigoplus}\e
\br(\be X_{K_v}\be)\! @>>>\ng\underset v\to{\bigoplus}\e
H^1(\gv,P)\! @>>>\e 0\endCD $$
where the direct sums extend over all primes of $K$, the vertical
maps are the natural ones and $g=\oplus_v g_v$. This diagram is of the
form considered at the end of Section 2, and we may therefore
apply to it Proposition 1.6 above. Before doing so, however, we call upon

\proclaim{Lemma 2.2}{\rm {(a)}} There is an exact sequence
$$0\ra\br(K)\overset\eta\to{\longrightarrow}\bigoplus_v\br(K_v)
\overset\sum\roman{inv}_v\to{\longrightarrow}\Bbb Q/\Bbb Z\ra 0.$$
{\rm {(b)}} There is an exact sequence
$$0\ra\br(X)\ra\br(X_K)\ra\bigoplus_{v\notin S}\br(X_{K_v})\,,$$
where $S$ is the set of primes of $K$ not corresponding to a point
of $U$.
\endproclaim
{\it Proof}. Assertion (a) is one of the major theorems of class
field theory. See [47, \S 11]. Assertion (b) is proved in [28,
Lemma 2.6].\qed

We now apply Proposition 1.6 to the diagram above using the preceding
lemma. We get an exact sequence
$$0\ra P(K)/\pic(X_K)\ra\krn{\bar{\eta}}\ra\br(X)'\ra
{\cyr W}(P)\ra\cok{\bar{\eta}},$$
where $\bar{\eta}\:\br(K)\ra\bigoplus_v\br(K_v)\big/\,\img g\;$ is
induced by $\eta$ and
$${\br}(X)'={\krn}\!\Biggl[@,@,\br(X)\rightarrow
\bigoplus_{v\in S}\br(X_{K_v})@,@,\Biggr ].$$
Now by Lemma 1.1, the order of $P(K)/\pic(X_K)$ equals
$(\x/\y)\cdot\lb A(K):\pic^0(X_K)\rb$. Regarding the kernel and
cokernel of $\bar{\eta}$, the following holds.

\proclaim{Proposition 2.3} We have
$$\lb\krn\bar{\eta}\rb =\tprod\e\xv/\te\,,$$
and the map  $\sum\roman{inv}_v\,\:
\bigoplus\br(K_v)\ra\Bbb Q/\Bbb Z$ induces an isomorphism
$$\cok{\bar{\eta}}\,\simeq\,\Bbb Q\e/\e\te^{\be -1}\Bbb Z\,.$$
\endproclaim
{\it Proof}. By combining Lemmas 1.7, 2.1 and 2.2 we obtain an
exact sequence
$$0\ra\krn{\bar{\eta}}\ra\bigoplus\br(K_v)_{\xv\be -\roman{tor}}\overset
\sum\roman{inv}_v\to{\longrightarrow}\Bbb Q/\Bbb Z
\ra\cok{\bar{\eta}}\ra 0.$$
Now for each $v$ the invariant map
$\roman{inv}_v$ induces an isomorphism
$\br(K_v)_{\xv\be -\roman{tor}}\simeq\xv^{-1}\e\Bbb Z/\Bbb Z$,
and it follows that the kernel and cokernel
of the middle map in the above exact sequence may be identified with
the kernel and cokernel of the map $\Sigma$ considered in Section 2. The
proposition now follows from Lemma 1.2.\qed
\medskip

We summarize the results obtained so far.

\proclaim{Corollary 2.4} There is an exact sequence
$$0\ra T_0\ra T_1\ra\br(X)'\ra{\cyr W}(P)\ra\Bbb Q/\te^{\be -1}\Bbb Z\,,$$
where $T_0$ and $T_1$ are finite groups of orders
$$\align \lb T_{0}\rb &=(\x/\y)\cdot\lb A(K):\pic^0(X_K)\rb\\
\lb T_{1}\rb &= \tprod\e\xv/\te\,.\endalign$$
\endproclaim

We now prove
\proclaim{Theorem 2.5} Assume that the integers $\yv$ are
relatively prime in pairs. Then
$$A(K)=\pic^0(X_K).$$
\endproclaim
{\it Proof}. There is an exact commutative diagram
$$\CD
0 @>>>A(K)/\pic^0(X_K) @>>>\br(K)\\
      @.     @VVV           @VV\eta V\\
      0 @>>>\underset{\text{all $v$}}\to{\bigoplus}
A(K_v)/\pic^0(X_{K_v}) @>>>
\underset{\text{all $v$}}\to{\bigoplus}\br(K_v)\,,\endCD $$
in which the vertical maps are the natural ones and the nontrivial
horizontal maps are induced by the maps $P(K)/\pic(X_K)\ra\br(K)$ and
$P(K_v)/\pic(X_{K_v})\ra\br(K_v)$ coming from (5) and (6).
Arguing as in the proof of Lemma 2.1
(using Lemma 1.4 in place of Lemma 1.5 there), we conclude that
for each $v$ the image of $A(K_v)/\pic^0(X_{K_v})$ in $\br(K_v)$ equals
$\br(K_v)_{\yv\be -\roman{tor}}$. It follows that $A(K)/
\pic^0(X_K)$ injects into $\krn\bar{\eta}^{\e\prime}$, where

$$\bar{\eta}^{\e\prime}\,\:\br(K)\ra\bigoplus_v\br(X_{K_v})/\e\br(K_v)_
{\yv\be -\roman{tor}}$$
is induced by $\eta$. Now arguing as in the proof of Proposition 2.3, we
see that the order of
$\krn\bar{\eta}'$ equals $\prod\e\yv/\ue$. By hypothesis this number is
1, whence the theorem follows.\qed

We now recall from the Introduction the integer $d$, which was defined
to be the number of primes $v$ of $K$ for which
$\xv=2\e\yv$. As noted
just before the statement of Lemma 1.2, $\te/\ue=\text{1 or 2}$
according as $d=0$ or $d\geq 1$. Now we observe that

$$\align \frac{\tprod\e\xv}{\te} & =\fracwithdelims(){\te}{\ue}^{-1}
\prod\e(\xv/\yv)\cdot\frac{\tprod\e\yv}{\ue}\\
& =2^e\,\frac{\tprod\e\yv}{\ue}\,,\endalign$$
where
$$e={\roman{max}}(0, d-1).$$
Consequently if the integers $\yv$ are relatively prime in pairs,
then $\prod\e\xv/\te = 2^e$. Thus Corollary
2.4 and Theorem 2.5 together imply 

\proclaim{Corollary 2.6} Assume that the integers $\yv$ are
relatively prime in pairs. Then there is an exact sequence
$$0\ra T_0\ra T_1\ra\br(X)'\ra{\cyr W}(P)\overset\phi\to
{\longrightarrow}\Bbb Q/\te^{\be -1}\Bbb Z\,,$$
in which $T_0$ and $T_1$ are finite groups of orders
$$\align \lb T_{0}\rb &=\x/\y\\
\lb T_{1}\rb &= 2^e\e ,\endalign$$
where $e={\roman{max}}(0, d-1)$.
\endproclaim

We turn now to the problem of relating ${\cyr W}(P)$ to ${\cyr W}(A)$.
There is an exact commutative diagram
$$\minCDarrowwidth{.5cm}
\CD 0 @>>>\Bbb Z/\y\e\Bbb Z @>>> H^1(\g,A) @>>> H^1(\g,P) @>>> 0\\
      @.     @VV D V           @VVV        @VVV            @.\\
      0 @>>>\underset{\text{all $v$}}\to{\bigoplus}\e\Bbb Z/\yv\e
\Bbb Z @>>>\underset{\text{all $v$}}\to{\bigoplus} H^1(\gv,A) @>>>
\underset{\text{all $v$}}\to{\bigoplus} H^1(\gv,P) @>>> 0\,,\endCD $$
in which $D$ is the diagonal map of Section 2 and the rows come from the
cohomology sequences of
$$0\ra A\ra P\overset\dgr\to{\longrightarrow}\Bbb Z\ra 0$$
over $K$ and over $K_v$. Applying the snake lemma to the above diagram
yields the exact sequence
$$0\ra\krn D\ra {\cyr W}(A)\ra{\cyr W}(P)\ra\cok D.\tag 7$$
Let $T_2=\krn D$. Then, by Lemma 1.3,
$$\lb T_2\rb=\y/\ue\,.\tag 8$$
Further, the order of $\cok D$ is $\prod\yv/\ue$. Consequently the
following holds

\proclaim{Proposition 2.7} Suppose that the integers $\yv$ are relatively
primes in pairs. Then there is an exact sequence
$$0\ra T_2\ra {\cyr W}(A)\overset\rho\to{\longrightarrow}{\cyr W}(P)
\ra 0\,,$$
in which $T_2$ is a finite group of order
$$\lb T_{2}\rb =\y/{\tsize{\prod}}\e\yv.$$
\endproclaim

In what follows we will view $T_2$ as a subgroup of ${\cyr W}(A)$ by
identifying it with its image in ${\cyr W}(A)$ under the map in
Proposition 2.7.

Putting together Corollary 2.6 and Proposition 2.7, we obtain

\proclaim{Corollary 2.8} Assume that the integers $\yv$ are
relatively prime in pairs. Then there is an exact sequence
$$0\ra T_0\ra T_1\ra\br(X)'\ra{\cyr W}(A)/T_2\ra\Bbb Q/\te^{\be -1}
\Bbb Z,$$
in which $T_0$, $T_1$ and $T_2$ are finite groups of orders
$$\align \lb T_{0}\rb &=\x/\y\\
\lb T_{1}\rb &= 2^e\\
\lb T_2\rb &=\y/{\tsize{\prod}}\e\yv\,,\endalign$$
where $e={\roman{max}}(0, d-1)$.
\endproclaim

It remains only to compute $T_3$, the image of ${\cyr W}(A)/T_2$ in
$\Bbb Q/\te^{\be -1}\Bbb Z$ under the map in Corollary 2.8, or
equivalently, the image of the composite map ${\cyr W}(A)\overset\rho
\to{\longrightarrow}{\cyr W}(P)\overset\phi\to{\longrightarrow}
\Bbb Q/\te^{\be -1}\Bbb Z$,
where $\phi$ and $\rho$ are the maps in Corollary 2.6 and Proposition
2.7, respectively. We will show that $T_3$ has the order stated in
the Introduction by generalizing [28, Lemma 2.11].

We begin by recalling from [28, Remark 2.9] the explicit description of
the map $\phi\,\:{\cyr W}(P)\ra\Bbb Q/\te^{\be -1}\Bbb Z$.

We write $E$ for the canonical map $\iv(X_{\kb})\ra\pic(X_{\kb})=P(\kb)$.
Represent $\alpha\in{\cyr W}(P)$ by a cocycle $a\in Z^1(\g,P(\kb))$, and
let
$\frak a\in C^1(\g,\iv(X_{\kb}))$ be such that $E(\frak a)=a$. Then
$\partial\e\frak a\in Z^2(\g,\krn E)=Z^2(\g,R(X_{\kb})^*/\kb^*)$ and,
because $H^3(\g,\kb^*)=0$, it can be pulled back to
an element $f\in Z^2(\g,R(X_{\kb})^*)$ (here $\partial=\text{boundary map}
$). On the other hand, $\roman{res}_v
(a)=\partial\e a_v$ with $a_v\in C^{\e 0}(\gv,P(\kb_v))$. Let $\frak a_
{\e v}\in
C^{\e 0}(\gv,\iv(X_{\kb_v}))$ be such that $E(\frak a_{\e v})=a_v$. Then
$\roman{res}_v(\frak a)=\partial\e\frak a_{\e v} + \roman{div}(f_v)$
with $f_v\in C ^1(\gv,R(X_{\kb_v})^*)$, and $\roman{res}_v f/\partial
f_v\in Z^2(\gv, \kb_v^*)$. Let $\gamma_v$ be the class of
$\roman{res}_v f/\partial f_v$ in $\br(K_v)$. Then
$$\phi(\alpha)=q\Biggl(\sum_v\roman{inv}_v(\gamma_v)\Biggr)\,,$$
where $q$ is the canonical map $\Bbb Q/\Bbb Z\ra\Bbb Q/\te^{\be -1}
\Bbb Z$ induced by the identity map on $\Bbb Q$.

We note that, for any divisor $\frak c_{\e v}$ on $X_{K_v}$
such that neither $\roman{res}_v f$ nor $\partial f_v$ has a zero or
a pole in the support of $\frak c_{\e v}$,
$$(\roman{res}_v f)(\frak c_{\e v})/\partial f_v(\frak c_{\e v})=(\dgr\,
\frak c_{\e v})
\,\roman{res}_v f/\partial f_v$$
(see for example [21, \S 4]). Here $(\roman{res}_v f)(\frak c_{\e v})=
f(\frak c_{\e v})\in
Z^2(\gv,\kb_v^*)$ is the value at $(\roman{res}_v f, \frak  c_v)$ of
the cup-product pairing
$$Z^2(\gv,R(X_{\kb_v})^*)\times Z^{\e 0}(\gv,\iv(X_{\kb}))\ra Z^2(\gv,\kb_v^*)$$
which is induced by the evaluation pairing $R(X_{\kb_v})^*\times
\iv(X_{\kb_v})\ra \kb_v^*$, and similarly for $\partial f_v(\frak c_{\e
v})$.
Since $\partial f_v(\frak c_{\e v})=\partial(f_v(\frak c_{\e v}))$ with
$f_v(\frak c_{\e v})\in C^1(\gv,\kb_v^*)$, we see that
$(\dgr\,\frak c_{\e v})\,\gamma_v$ is represented by $f(\frak c_{\e v})$.
Now choose a divisor $\frak c_{\e v}$ on $X_{K_v}$ of degree $\xv$ such that
neither $f$ nor $\partial f_v$ has a zero or a pole in the
support of $\frak c_{\e v}$ (this is always possible; see [19, App. 2]).
Then $\xv\gamma_v$ is represented by $f(\frak c_{\e v})$.

We now recall from [29, Remark I.6.12] the definition of a pairing
$$<\;,\;>\,\:{\cyr W}(A)\times{\cyr W}(A)\ra\Bbb Q/\Bbb Z\tag 9$$
which annihilates only the divisible part of ${\cyr W}(A)$.

Let $\alpha\in{\cyr W}(A)$ be represented by $a\in Z^1(\g,A(\kb))$,
and let $\roman{res}_v(a)=\partial\e a_v$ with $a_v\in Z^{\e 0}
(\gv,A(\kb_v))$.
Write
$$\align a & =E(\frak a),\qquad\frak a\in C^1(\g,\iv^0(X_{\kb}))\\
a_v & =E(\frak a_{\e v}),\qquad\frak a_{\e v}\in C^{\e 0}(\gv,\iv^0
(X_{\kb_v})).\endalign$$
Then $\roman{res}_v(\frak a)=\partial\e\frak a_{\e v} + \roman{div}
(f_v)$ in $C^1(\gv,\iv^0(X_{\kb_v}))$ with $f_v\in C ^1(\gv,R(X_{
\kb_v})^*)$. Moreover, $\partial\e\frak a=\roman{div}(f)$ with
$f\in Z^2(\g,R(X_{\kb})^*)$. Let $\alpha'$ be a secod element of
${\cyr W}(A)$ and define $\frak a',\frak a_{\e v}',f_v',f'$ as for
$\alpha$. Set
$$<\alpha, \alpha'>=\sum_v\roman{inv}_v(\gamma_v),\qquad
\gamma_v=\text{class of ($f_v'\cup\roman{res}_v(\frak a)-\frak a'_v
\cup\roman{res}_v f\e $)},$$
where $\cup$ denotes the cup-product pairing induced by the
evaluation pairing mentioned earlier. This definition is independent
of the choices made, and $f,\frak a'_v, f_v'$ and $\frak a$ can be chosen
so that $f_v'(\roman{res}_v(\frak a))$ and $f(\frak a'_v)$
are defined. Moreover, $<\alpha,\alpha>=0$.

Now let $\varepsilon\,\:\Bbb Q/\Bbb Z\overset\sim\to{\longrightarrow}
\Bbb Q/\te^{\be -1}\Bbb Z$
be the natural isomorphism which is induced by multiplication by
$\te^{\be -1}$ on $\Bbb Q$.

\proclaim{Proposition 2.9} Let $\beta$ be a generator of $T_2\subset
{\cyr W}(A)$ and let $\beta'=(\te/\ue)\beta$. Then the composite
$${\cyr W}(A)\overset\rho\to{\longrightarrow}
{\cyr W}(P)\overset\phi\to{\longrightarrow}\Bbb Q/\te^{\be -1}\Bbb Z$$
is $\alpha\mapsto\varepsilon(<\alpha,\beta'>)$.
\endproclaim
{\it Proof}. Let $\alpha\in{\cyr W}(A)$ and define $\frak a,\frak a_{\e v},
f_v$ and $f$ as above. Then $\phi(\rho(\alpha))=\sum\roman{inv}_v(
\gamma_v)$ where $\xv\gamma_v$ is represented by $f(\frak c_{\e v})$ for
some divisor $\frak c_{\e v}$ of degree $\xv$ on $X_{K_v}$.

On the other hand, we can choose $\beta'\in{\cyr W}(A)$ to be represented
by $b'=E(\frak b')$ where $\frak b'=\te\e\partial\e P$, $P$ any closed
point on $X_{\kb}$. Define
$$\frak b'_v=\te\e P-(\te/\xv)\frak c_{\e v}\in C^{\e 0}(\gv,
\iv^0(X_{\kb_v})).$$
Then $b'_v=E(\frak b'_v)$ satisfies
$$\partial\e b'_v= E(\te\e\partial\e P)=E(\roman{res}_v\frak b')=
\roman{res}_v(b').$$
Further, $\roman{res}_v(\frak b')=\partial\e\frak b'_v$, which means
that we may choose $f'_v=1$. Thus $<\alpha ,\beta'>=-\sum\roman{inv}_v(
\gamma_v^{\e\prime})$, where $\gamma_v^{\e\prime}$ is represented by
$f(\frak b_v')=f(\te\e P-(\te/\xv)\frak c_{\e v})=f(P)^{\te}/
\gamma_v^{\e\te}$. Consequently
$$\align\varepsilon(<\alpha,\beta'>)& =\varepsilon\biggl(\te\sum_v
\roman{inv}_v(\gamma_v)\biggr)\\
& =q\biggl(\sum_v\roman{inv}_v(\gamma_v)\biggr)=\phi(\rho(\alpha))\,,
\endalign$$
as claimed.\qed

\proclaim{Corollary 2.10} Assume that ${\cyr W}(A)$ contains no
nonzero infinitely divisible elements. Then
$$\lb T_3 \rb =\frac{\y/\ue}{\bigl(\frac{\te}{\ue},\frac{\y}{\ue}
\bigr)}.$$
\endproclaim
{\it Proof}. Under our hypothesis the pairing (9) is nondegenerate,
and the proposition shows that $T_3$ is isomorphic to the dual
of $<\beta'>=(\te/\ue)\e T_2$. The corollary now follows from (8).\qed

The Main Theorem of the Introduction may now be obtained by
combining Corollaries 2.8 and 2.10.

\heading 3. Concluding Remarks\endheading

As shown in the proof of Theorem 2.5, the order of
$A(K)/\pic^0(X_K)$ divides $\prod\e\yv/\ue$. We believe that
it is always equal to $\prod\e\yv/\ue$. This is equivalent to the
following plausible statement. Consider the nondegenerate pairing
$$\bigoplus_v\Bbb Z/\yv\e\Bbb Z\times\bigoplus_v A(K_v)/\pic^0(X_{K_v})
\ra\Bbb Q/\Bbb Z\,,$$
which is the sum of the local pairings
$$\Bbb Z/\yv\e\Bbb Z\times A(K_v)/\pic^0(X_{K_v})\ra\Bbb Q/\Bbb Z$$
induced by Lichtenbaum duality (see the discussion preceding the
statement of Lemma 1.4). Then, relative to this pairing, the
image of the diagonal map
$$D\,\:\Bbb Z/\y\e\Bbb Z\ra\bigoplus_v\e\Bbb Z/\yv\e\Bbb Z$$
is the exact annihilator of the image of the diagonal map
$$\Cal D\,\: A(K)/\pic^0(X_{K})\hookrightarrow\bigoplus A(K_v)/
\pic^0(X_{K_v}).$$

Assuming (for simplicity) that ${\cyr W}(A)$ is finite, the above
conjectural statement implies in addition that the exact sequence

$$0\ra\krn D\ra {\cyr W}(A)\ra{\cyr W}(P)\ra\cok D\ra\cok\lambda\ra
\cok\mu\ra 0$$
(which is the continuation of the exact sequence (7). Here $\lambda$
and $\mu$ are the natural maps $H^1(\g,A)\ra\bigoplus H^1(\gv,A)$
and $H^1(\g,P)\ra\bigoplus H^1(\gv,P)$) splits into two short exact
sequences
$$0\ra\krn D\ra {\cyr W}(A)\ra{\cyr W}(P)\ra 0$$
and
$$0\ra\cok D\ra\cok\lambda\ra\cok\mu\ra 0,$$
the second of which is the dual of
$$0\ra\pic^0(X_K)\ra A(K)\ra\img\Cal D\ra 0.$$

Thus if the above conjecture is true, then we can give a complete answer to
Grothendieck's question stated in the Introduction.

\enddocument